\documentclass[11pt]{report} 

\usepackage{amsfonts,
            amssymb, young}

\def \dfll {\leaders \hbox to 1em {\hss.\hss}\hfill}
\def\beq {\begin{equation}}
\def\eeq {\end{equation}}
\newtheorem{theorem}{Theorem}
\newtheorem{propos}[theorem]{Proposition}

\newtheorem{qes}{Question}
\newtheorem{lemma}{Lemma}

\def\byg{\begin{Young}}
\def\eyg{\end{Young}}
\def\bsyg{\begin{SmallYoung}}
\def\esyg{\end{SmallYoung}}
\def\btid{\begin{Tabloid}}
\def\etid{\end{Tabloid}}
\def\bstid{\begin{SmallTabloid}}
\def\estid{\end{SmallTabloid}}

\def\blm{\begin{lemma}}
\def\elm{\end{lemma}}
\def\bdf{\begin{Defi}}
\def\edf{\end{Defi}}
\def\btm{\begin{theorem}}
\def\etm{\end{theorem}}
\def\bpp{\begin{propos}}
\def\epp{\end{propos}}
\def\bQ {\begin{qes}}
\def\eQ {\end{qes}}
\def\btm{\begin{theorem}}
\def\etm{\end{theorem}}
\def\ben{\begin{enumerate}}
\def\een{\end{enumerate}}

\def\ep{\ \hfill{\rule {2.5mm}{2.5mm}}\smallskip}

 \newcommand{\sext}[1] 
{\mbox{\raisebox{.2ex}{$\scriptstyle \bigwedge^{\!#1}$}}\mkern-1mu}

\def \isto {\widetilde{\longrightarrow}}

\def \Cob {\mbox{\boldmath${\cal C}ob$\unboldmath}}

\def \sl {{\mathfrak s}{\mathfrak l}}

\newcommand{\picbox}[1]{\begin{equation}
\epsfbox{#1.eps}\nonumber\end{equation}}

\newcommand{\TO}[2]{\stackrel {\mbox{#1}}{\hbox to #2pt{\rightarrowfill}}}

\def\thrafill{$\mathsurround=0pt \mathord- \mkern-6mu 
\cleaders\hbox{$\mkern-2mu
\mathord- \mkern-2mu$}\hfill \mkern-6mu\mathord\twoheadrightarrow$}
\newcommand {\onto} [1]{\hbox to #1pt{\thrafill}}

\usepackage{amsfonts,
            amssymb}
\usepackage{
            xypic}      

\input{epsf.sty}    
\def\opc{{\scriptstyle =}}

\def\drulefill{$\mathord\opc\mkern-3mu\cleaders
                  \hbox{$\mkern-2mu\mathord\opc\mkern-2mu$}\hfill\mkern-2mu \mathord\opc$}

\def\sdwa{$\raise+5pt\hbox{$\scriptstyle\rangle\mkern-3.5mu\rangle$}
\mkern-9.55mu\raise-0pt\hbox{$\scriptstyle\langle\mkern-3.5mu\langle$}
\mkern-6.95mu\raise-5pt\hbox{$\scriptstyle\rangle\mkern-3.5mu\rangle$}
\mkern-11.9mu\raise-11.1pt\hbox{$\scriptstyle\bigtriangledown$}
$}

\def\dov#1{\vbox{\ialign{##\crcr\noalign
             {\kern-3pt\nointerlineskip}\drulefill \crcr\noalign
             {\kern0pt\nointerlineskip}
             $\hfil\displaystyle{#1}\hfil$\crcr}}}

\def\zer#1{\vbox{\ialign{##\crcr\noalign
             {\kern-3pt\nointerlineskip}
             {$\hfil\;\,\scriptscriptstyle\oslash\hfil\!$} \crcr\noalign
             {\kern.5pt\nointerlineskip}
             $\hfil\displaystyle{#1}\hfil$\crcr}}\!}

\def\ste#1{\vbox{\ialign{##\crcr\noalign
             {\kern-3pt\nointerlineskip}
             {$\hfil\;\,\scriptscriptstyle\bowtie\hfil\!$} \crcr\noalign
             {\kern.5pt\nointerlineskip}
             $\hfil\displaystyle{#1}\hfil$\crcr}}\!}

\newcommand{\Z}{{\mathbb Z}}

\newcommand{\R}{{\mathbb R}}
\newcommand{\Cc}{{\mathbb C}}
\newcommand{\N}{{\mathbb N}}
\newcommand{\Q}{{\mathbb Q}}

\newcommand{\Pp}{{\mathbb P}}
\newcommand{\FF}{{\mathbb F}}
\newcommand{\I}{{\mathbb I}}

\newcommand{\posit}[2]
{\raise -1.4ex\hbox{${\textstyle #1}\atop {\stackrel{\uparrow}{#2}}$}}

\newcommand{\vs}[1]{\vspace*{#1mm}}

\def \lz  {\langle}
\def \rz  {\rangle}

\def\DP{\displaystyle}

\newcommand{\ext}[1] {\mbox{\raisebox{.4ex}{$\bigwedge^{\!#1}$}}\mkern-1mu}

\def\Ii {\mbox{\raise .4 ex\hbox{$\int$}$\!\! I$}}

\begin{document}

\begin{center}
\ 

\vs{6}

{\bf\LARGE $p$-Modular TQFT's and Torsion} 

\vs{13}

{\large \sc Thomas Kerler}
\vs{9}

 { \small \sc The Ohio State University \\ Columbus, OH, USA}

\vs{9}

 September 6th, 2001

\vs{4} 

{\small presented  in Seminar/Workshop on }

\vs{6} 

{\bf \large ``Invariants of knots and 3-manifolds" }

\vs{6}

{\small at the } 

\vs{8}

 {  \sc Research Institute for Mathematical Sciences}

\vs{6}

{   \sc Kyoto, Japan}

\end{center}

\newpage
 
\section*{Non-Perturbative Quantum Invariants:}

Let $p\geq 5$ be a prime, \ $\zeta_p$ primitive $p$-th root of unity.  

${\cal C}_p$ semisimple category over $\Z[\zeta_p]$ given by

\quad even ($SO(3)$) part of $\overline{U_{\zeta_p}({\mathfrak sl}_2)-mod}$ 
\medskip

\subsection*{Quantum-Invariants:}  
\ \ $M$ compact, oriented 3-manifold with 
$\partial M = \emptyset$. 
\bigskip 

\noindent
{\footnotesize\bf [Reshetikhin, Turaev]}\ \   ${\cal V}_{\zeta_p}^{RT}(M)\in\Cc\,$\  
w/ \  ${\cal V}_p^{RT}(S^3)=1$.
\bigskip 

\noindent
{\footnotesize\bf  [Murakami, Masbaum, Roberts]}\ \ \  ${\cal V}_{\zeta_p}^{RT}(M)\in \Z[\zeta_p]$.
\bigskip 

\noindent
{\small\bf [Ohtsuki]}\ \ \  For $\beta_1(M)=0$ have $\lambda_j^O(M)\in
\Z[\frac 12,\ldots, \frac 1{2j+1}]$
$$
{\cal V}_{\zeta_p}^{RT}(M)=\sum_{j=0}^{\frac {p-3}2}\lambda_j^O(M)(\zeta_p-1)^j
\;+\;{\cal O}((\zeta_p-1)^{\frac {p-1}2})\;.
$$

\noindent
{\small\bf [Murakami]}\qquad $\lambda_{Casson}=6\lambda_1^O$
\smallskip

\subsection*{TQFT:} 

{\small\bf [Atiyah]} TQFT \ = \ Functor ${\cal V}:\Cob\,
\longrightarrow\,{\sf R}-mod$

\bigskip

\noindent
{\footnotesize\bf [Reshetikhin, Turaev]}\ \ \ 
Constructions ${\cal V}_{\zeta_p}^{RT}$ for ${\sf R}=\Cc$. 

\bigskip

\noindent
{\small\bf [Gilmer]} \ \ \ For ${\sf R}=\Z[\zeta_p]$ for ``grounded'' cobordisms. 

\bigskip

\noindent
{\bf [K]}  \ \ \ For ${\sf R}=\Z[\zeta_5]$ ``half-projective'' TQFT. Basis.

\bigskip

\newpage

\section*{Topological Quantum Field Theory:}

\paragraph{Cobordisms:} $\partial M=\partial^+M\sqcup\partial^-M$ w/ homeom.
\smallskip

\ \ $\psi^{\pm}:\,\partial^{\pm}M\,\isto\,\Sigma_{\pm}$ to model surfaces.
\smallskip

Homeom. class $[M,\psi^{\pm}]:\,\Sigma_-\longrightarrow \Sigma_+$ 
(+Framings).  
\picbox{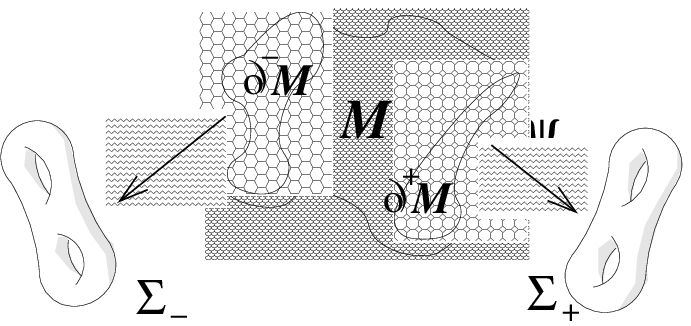}
\paragraph{Category $\Cob$:} Composition by ``Gluing'' over surfaces. 
 
\paragraph{TQFT Functor:} Assigns ${\cal V}(\Sigma)\in {\sf R}-mod$ and for

 $M:\Sigma_-\to\Sigma_+$\ \  have\  \  ${\cal V}(M)\;\in\;
Hom_{\sf R}({\cal V}(\Sigma_-), {\cal V}(\Sigma_+))$. 
\vs{-1}


\noindent 
 {\bf [K]}``Half-Projective'': \ \ \ \ 
 ${\cal V}(N\circ M)={\sf x}^{\mu(N,M)}{\cal V}(N){\cal V}(M)$
\vs{-1}

\noindent 
\ \ with  \ 
$\mu(N,M)=rk(H_1(N\circ M)\stackrel{\delta}{\longrightarrow}
H_0(M\cap N))$, \  ${\sf x}\in{\sf R}$
\vs{-1}

\ \ \ \ Ordinary  Composition:  \ \ ${\sf x}\;=\;1$. 

\paragraph{Mapping Class Group Representation:} \ \ 
\begin{eqnarray}
\Gamma_g=\pi_0(\mbox{\it \small Homeo}^+(\Sigma))\;\;&\stackrel{\cong}
{\longrightarrow}&\;\;Aut(\Sigma)\;\;\stackrel{{\cal V}}{- - - \to}GL({\cal V}(\Sigma_g))\nonumber \\
\; [\psi] \;\;\;&\mapsto&\;\Bigl [\Sigma\times[0,1],\psi\sqcup id\Bigr]
\nonumber 
\end{eqnarray}

\newpage

\section*{Main Question:}

{\large \bf What is the  TQFT content of the $\lambda_j^O$ ?}
\smallskip

\noindent 
Expansions \ \     $\DP {\cal V}_{\zeta_p}^{RT}(M)\;=\;\sum_j
{\cal V}_{p,j}^{RT}(M){\sf y}^j$\ \  for \ \ $M\in\Cob$

via \ \ \ \ $\Z[\zeta_p]\,\;\;\onto{4.5}\;\;\quad\FF_p[\zeta_p]\;\cong \;\FF_p[{\sf y}]\Big/ {\sf y}^{p-1}$
\medskip 

w/ \ \ \  $\FF_p=\Z/p\Z$ \ \ and \ \ ${\sf y}=\zeta_p-1$

\subsection*{Applications \& Remarks:}

\begin{itemize}

\item Computations via TQFT-techniques

(Heegaard splits, Mapping tori, Seifert spaces, ...) 

\item Extensions of $\lambda_j^O$ to $\beta_1(M)\,>\,0$.

(Reidemeister Torsion vs. Casson Invariant)

\item Filtration  
${\cal M}_{p,k}\,=\,\DP \bigcap_{j=0}^{k}ker \Bigl({\cal V}_{p,j}^{RT}\Big |_{\Gamma_g}\Big)$ of {\em mapping

\ \  class group} $\Gamma_g$ \ \ 
$\stackrel{\mbox{\bf ?}}{\Longleftrightarrow}$ \ Other (f.t.) Filtrations.  

\item Dependence on $p$: \ $dim({\cal V}_{\zeta_p}^{RT}(\Sigma_g)) 
\sim p^{3(g-1)}\;\mbox{if}\,g>1
$

\ \ $\Longrightarrow$ \ \ \  NO fixed universal TQFT !

\item Even for FIXED $p$, operators ${\cal V}_{p,j}^{RT}(M)$ over $\FF_p$ are
 
NOT liftable to $\Z[\ldots]$ as TQFT's. 
\end{itemize}

\newpage

\section*{Some Results - Overview:}

\paragraph{General primes $p\geq 5$: }

\begin{itemize}

\item Construct Family of {\em irreducible} TQFT's over $\FF_p$

$\dov {\cal V}_p^{(k)}$ \ 
($\dov {\cal U}_p^{(k)}$) for $p$ odd prime and $0\,<\,k\,<\,p(-3)\,$. 

\item Resolutions  
$\,0 \leftarrow \dov {\cal V}_p^{(k)}\leftarrow   {\cal V}_p^{(k)}
\leftarrow    {\cal V}_p^{(2p-k)}\leftarrow \ldots \,$ 

into $\Z$-liftable TQFT's \ \  ${\cal V}_p^{(j)}\leftarrow \!\!\!\!
\leftarrow {\cal V}_{\Z}^{(j)}$. 

\item Image of Alexander Polynomial $\dov \Delta_{p,\chi}(M)\in\FF_p[\zeta_p]$ 

 from traces over $\dov {\cal V}_p^{(k)}$ of covering cobordism. 

\item Similar results for $p$-modular  $S_n$ representations. 

\end{itemize}

\paragraph{Fibonacci Case $p = 5$: }

\begin{itemize}

\item ${\cal V}^{RT}_{\zeta_5}\,$ =\ half-projective, with \ 
 ${\sf x}=\zeta_5-\zeta_5^{-1}\in \Z[\zeta_5]$. 

\item \ \ Identification: \ \   \ \   $\DP {\cal V}^{RT}_{5,0}\;\;=\;\;\dov {\cal U}_5^{(1)}$

\item \ \ $2\dov \Delta_{5  ,\chi}(M)\;=\;det(ker(\chi))\,+\,3{\sf y}{\cal V}_{\zeta_5}^{RT}(M)
\;\,+\;\,{\cal O}({\sf y}^3)$

\item \ \ \ \ \ \ $cut(M)\;\leq\;{\mathfrak o}_5(M)$\ \ \ (joint w/ Gilmer)

\item \ \  \ \  $\lambda_{Casson}(M_{D(\psi({\cal C}_k))})=\lz \Omega_{\Lambda},\,
\Pp^{\psi}_k \, \Omega_{\Lambda}\rz\;\;\,{\rm mod}\, 5$

\ \ \  $ \Omega_{\Lambda}\in\ext g H_1(\Sigma_g)$, 
\ ${\cal C}_k$ $k$-th bounding curve; $\psi\in\Gamma_g$; 
$D({\cal C})$ Dehn twist @ ${\cal C}$;
$M_{\gamma}$ Heegaard glue; $\Pp^2_k=\Pp_k$.  

\item $\lambda_j^O \,{\rm mod}\, 5$ are finite type, ${\cal V}^{RT}_{0,5}=det(H_1(M))$, ... 
\end{itemize}

\section*{The Alexander Polynomial:}

\noindent
\parbox[t]{2.6in}
{\vs{-3}

\noindent 
$M$ = compact, oriented,  $\beta_1(M)>0$
\vs{3}

$\chi:H_1(M)\,-\!\!\!\!\!\!\to\!\!\!\!\to \,\Z$\ \   epimorphism}
\parbox[t]{1.8in}
{\vs{-9}

\ \ \ \ \ $
\begin{array}{c}
\widetilde M_{\chi}\\
\downarrow\\
M\\
\end{array}$

 cyclic covering. } 
\vs{2}

$H_1(\widetilde M_{\chi},\Z)\,$ is 
$\,\Z[t,t^{-1}]$-module:\ \  $t$ by  Deck transf. 
\vs{3}

\noindent
{\bf Definition:} \ \ $\Delta_{\chi}(M)\;\in\;\Z[t,t^{-1}]$ \ \ as follows
\begin{itemize}
\item  Maximal such that 
$E_1(H_1(\widetilde M_{\chi}))\subseteq \Delta_{\chi}(M)\Z[t,t^{-1}]$
\item $det(A(t))$ if $A(t)$ is $n\times n$-present. matrix of 
$H_1(\widetilde M_{\chi})$.
\item Choose symmetric (by P.D.) $\Longrightarrow$  Unique up to sign. 
\end{itemize} 
\vs{1}

\noindent
{\bf Thm \small[Milnor, Turaev]} Let $r_{\chi}(M)$ {\em Reidemeister Torsion}: 
\vs{-3}

 If $\partial M\neq \emptyset$ then 
$\displaystyle r_{\chi}(M)\,=\, {(t-1)}^{-1}\Delta_{\chi}(M)$
\vs{2}

If $\partial M = \emptyset$ then 
$\displaystyle r_{\chi}(M)\,=\, {(t-1)}^{-2}\Delta_{\chi}(M)$
\vs{2}

\noindent
{\bf Remarks:} \ \vs{-2}
\begin{itemize}
\item For  knot $K$:  $M_K=S^3-N(K)$ and $H_1(M_K)=\Z$. 
\item Multivariable $\Delta(M)\in \Z[H_1^{\it free}(M)]$\vs{-4}
\item $\Delta_{\chi}(1)=det(ker(\chi))$ \  \ w/ \  \  
       $det(H)$ \raise3pt\hbox{$=$} \raise6pt\hbox{$\left\{\begin{array}{cl} |H| & \mbox{if}\;<\infty\\
0 & \mbox{elsewise}\end{array}\right.$}
\item $\Delta_{\chi}=0$\ \ if\ \ $\chi:\pi_1(M)\,-\!\!\!\!\!\!\to\!\!\!\!\to  
\,\Z*\Z\,-\!\!\!\!\!\!\to\!\!\!\!\to \,\Z$. 
\item $\Delta_p(S^1\times\Sigma_g)=(t+t^{-1}-2)^g$ \ ...\ dependence on $\chi$!
\item {\bf Thm\ \footnotesize  [Meng, Taubes, Hutchings, Lee, Turaev, ...] }
\ \ \ {\em Seiberg-Witten Invariant} \ = \ (refined) {\em Torsion}. 
\end{itemize}

\section*{TQFT's from Jacobians:}

\ \ ${\cal V}^{FN}_{\Z}(\Sigma)\,=\,\raise-1pt\hbox{\em\Large H\,}^*
\Bigl(Hom(\pi_1(\Sigma), U(1))\Bigr)\;\cong\;\ext * H_1(\Sigma)$
\vs{2}

\noindent
{\bf Thm \small [Frohman, Nicas]:}  Extends to TQFT 
\vs{2}

\ \ \ \ ${\cal V}^{FN}_{\Z}\,:\;\Cob\;\longrightarrow\;\Z-mod$\ \ 

\begin{itemize}
\item $\Z/2\Z$-projectivity/framing-dependence.
\item Half-projective with ${\sf x}\;=\;0$. 
\item $\Gamma_g\,-\!\!\!\!\!\!\to\!\!\!\!\to \, Sp(2g,\Z)\,$ Representation.
\item \ \ \ ${\cal V}^{FN}_{\Z}(M)\;=\;det(H_1(M,\Z))\;$\ \ \ for \ \ \ 
$\partial M=\emptyset$. 
\end{itemize}

\noindent 
Let $\Sigma\subset M$ 2-sided, oriented surface, dual to $\chi\in H^1(M)$ 
\vs{0}

\noindent 
{\bf Thm \small [F,N]} Set $C_{\Sigma}=(M-N(\Sigma)):\,\Sigma\longrightarrow\Sigma$\ \  in\ \  $\Cob$.

\ \ $\DP \Delta_{\chi}(M)\,=\,\sum_{j=0}^{2g}
(-t)^{j-g}trace\Bigl({\cal V}^{FN}_{\Z}(C_{\Sigma})\Big|_{\sext j H_1(\Sigma_g)}
\Bigr)$
\vs{3}

\noindent 
{\bf Thm   [K]}\ \ \ ${\cal V}^{FN}_{\Z}\;\cong\;
{\cal V}^{Hennings}_{\cal N}$ \ \ \ w/ \ \  ${\cal N}=\Z/2\Z\ltimes \ext *\R^2$
\vs{2}

\noindent
  $\Longrightarrow$ \ {\bf [K]}\ Extension of Alexander-Conway Calculus to

\qquad\ \   3-manifolds presented by Surgery. 

\noindent 
{\bf Note:}
\begin{itemize}
\item {\bf\small [F,N]} \ \ TQFT's via $IH_*$ of  $PSU(n)$-Moduli Spaces. 
\item {\bf\small [Donaldson]}\ \  TQFT's via connections on non-trivial 
$SO(3)$-bundle (Casson-type) and   {\em SW}(voretx)-equations. 
\end{itemize}

\section*{Hard Lefschetz Decomposition:}

Jacobians are K\"ahler\ \   $\Rightarrow$\ \   
Inner product on $\ext *H_1(\Sigma)$.
\vs{1}

\noindent 
\ \ Symplectic 2-form $\omega\,=\,\sum_ia_i\wedge b_i\;\in\;  \ext 2H_1(\Sigma)$.
\vs{2}

\noindent 
{\bf $\sl(2,\Z)$-Action:} $E$, $F$, $H$  generators act on $\ext *H_1(\Sigma)$ as 
\vs{2}

$Hx=\;(deg(x)-g)x$, \ \ \ \ $Ex=x\wedge \omega$, \ \ \ \ $F=E^*$.
\vs{4} 

\noindent
{\bf Thm:}
\ \ \ \   
\raise-5pt
\hbox{$\sl(2,\Z)$ \ \ $\stackrel{\mbox{DUAL on ${\cal V}^{FN}_{\Z}(\Sigma_g)$}}
{\leftarrow
\!\!\!-\!\!\!-\!\!\!-\!\!\!-\!\!\!-\!\!\!-\!\!\!-\!\!\!-
\!\!\!-\!\!\!-\!\!\!-\!\!\!-\!\!\!-\!\!\!-\!\!\!-\!\!\!-
\!\!\!
\to}$\ \  $Sp(2g,\Z)$} 
\vs{3}

hence  \qquad \ $\DP \ext *H_1(\Sigma_g)\;=\;\bigoplus_{k=1}^{g} V_{k}\otimes W_{g,k}\;$
\vs{3}

\noindent 
where \qquad \ $V_k$ is $k-dim$ irred. $\sl_2$-module. 

\noindent
$W_{g,k}$ is irred. $Sp(2g,\Z)$-module; h.w.  
$\DP \varpi_{g,k}=\sum_{i=1}^{g-k+1}\epsilon_i$ 
\vs{3}

\noindent
{\bf Thm [K]:}\ \ \ \ \ \ $\DP  {\cal V}^{FN}_{\Z}\;\;\cong\;\;
\bigoplus_{k=1}^{\infty} V_{k} \otimes 
{\cal V}^{(k)}_{\Z}\;$

\begin{itemize}
\item ${\cal V}^{(k)}_{\Z}(\Sigma_g)\,=\,ker(F)\cap \ext {g-k+1}H_1(\Sigma_g)\,
\cong\, W_{g,k}$
\item ${\cal V}^{(k)}_{\Z}$ is irred. TQFT and ${\cal V}^{(k)}_{\Z}(\Sigma_g)=0$
for $k>g+1$. 
\item \ \ $\DP \Delta_{\chi}(M)\;\;=\;\;\sum_{k=1}^g \,[k]_{-t}\,  
trace({\cal V}^{(k)}_{\Z}(C_{\Sigma}))\;$ 
\vs{-2}

\qquad\qquad\qquad\qquad \qquad\qquad 
where   $\DP [n]_q\;=\;\frac {q^n-q^{-n}}{q-q^{-1}}\,$. 
\end{itemize}

\newpage 

\section*{$\FF_p=\Z/p\Z$-Reductions \& Reducibility:}

Let $p\geq 3$ be prime and ${\cal V}^{(j)}_p$ the $\FF_p$-reduction of 
${\cal V}^{(j)}_{\Z}$.
\vs{1.5}

\noindent
Inherits inner form 
$\,\lz\,,\rz_p:\, {\cal V}^{(j)}_p(\Sigma)\otimes {\cal V}^{(j)}_p(\Sigma)
\longrightarrow
\FF_p$
\vs{2.5}

\noindent
{\bf Note!} Forms $\lz\,,\rz_p$ degenerate $\Rightarrow$\ \ the  ${\cal V}^{(j)}_p$ are Reducible! 
\vs{-2.5}

\noindent
{\em Example:} 

\noindent
\ Let $v=E1-p(a_1\wedge b_1)=\omega - p(a_1\wedge b_1)\,\in\,\ext 2 H_1(\Sigma_p)$. 

\noindent
\ Have $Fv=p-p=0$ so that
$\,v\in  ker(F)\,=\,{\cal V}^{(p-1)}(\Sigma_p)$.

\noindent
\ For arbitrary $w\in {\cal V}^{(p-1)}(\Sigma_p)$ find
 
\noindent
\quad  $\lz v,w\rz_p=\lz E1,w\rz_p - \lz p(a_1\wedge b_1),w\rz_p =$

\qquad \  \ $=\lz 1,Fw\rz_p - p\lz (a_1\wedge b_1),w\rz_p\,=\,0$.

\vs{0}

\noindent
{\em Note also:} 

\ $\,E:\,{\cal V}^{(p+1)}_{\Z}(\Sigma_p)\;\to\;
{\cal V}^{(p-1)}_{\Z}(\Sigma_p)\,+\,p\ext 2 H_1(\Sigma_p,{\Z})$
\vs{2}

\noindent
\ \ so\ \ \ \ \  $\,E:\, {\cal V}^{(p+1)}_p(\Sigma_p)\;\to\;{\cal V}^{(p-1)}_p(\Sigma_p)$
\ \ \ \ \ well defined!
\vs{2}

\noindent
{\bf Thm [K]:} \ \ For $\,\,j\equiv k\,{\rm mod}\,p\,$ with $\,0< k<p\,$\  have 
\vs{2}

\ \ \ \quad $\DP \,E^k\,:\;{\cal V}^{(j)}_p\;\longrightarrow\;{\cal V}^{(j-2k)}_p$
\ \  \ \   well defined.
\vs{2}

All\ \ \  $E^k\neq 0$\ \ and \ \ $E^p\,=\,0$ in $p$-Reduction. 

\vs{1.5}
 
\noindent
$\Longrightarrow$ \ {\bf Complex:} \ \ For any $k$ with $0<k<p$ have 

\noindent 
${\cal C}(p,k)$\ \ \ $\;\ldots\,\longrightarrow\,{\cal V}^{(c_{i+1})}_p
\,\longrightarrow\,{\cal V}^{(c_i)}_p\,\longrightarrow\,\ldots\,$

$\ldots\,\to\, {\cal V}^{(2p+k)}_p
\,\longrightarrow\,{\cal V}^{(2p-k)}_p
\,\longrightarrow\,{\cal V}^{(k)}_p
\,\longrightarrow\,\dov {\cal V}^{(k)}_p\,\to\,0$
\vs{-1}

\noindent 
w/ $c_i=\scriptstyle \left\{\begin{array}{cl} ip+k & \mbox{$i$ even}\\
        (i+1)p-k & \mbox{$i$ odd}
\end{array}\right. $ \ \ \& \ \  maps $E^k$ or $E^{p-k}$. 
\vs{0.5}
 
\noindent
\ \ \ {\em Denote:} \ \ \ \ \ 
 $\dov {\cal V}^{(j)}_p$ the subquotient by $\lz \, ,\rz_p$-Null space.

\section*{Resolutions via Symmetric groups:}

{\bf Thm [K]}\ \ \ The sequences ${\cal C}(p,k)$ are {\em exact}. 
 \vs{1.5}

The $\dov {\cal V}^{(k)}_p$ are {\em irreducible} with 
{\em resolutions }  ${\cal C}(p,k)$. 
\vs{3}

\noindent
{\em Proof:} \ \  $Sp(2g,\Z)$-weight spaces 
${\cal W}^{(j)}(g,\lambda)\subset {\cal V}^{(j)}(\Sigma_g)$
\vs{1}

\ \ with  \ \ $\lambda\,=\,\sum_i\lambda_i\epsilon_i\,$ with $\lambda_j\in\{-1,0,1\}$
\vs{1}

Let $N(\lambda)=\{i:\lambda_i=0\}$ and $n(\lambda)\,=\,|N(\lambda)|$. 

\vs{1}

\noindent 
Via {\em Weyl Group}  $S_{n(\lambda)}$ acts on ${\cal W}^{(j)}(g,\lambda)$ by permut. 
$N(\lambda)$. 

\vs{2}

\noindent
{\bf Lem:}\ \ ${\cal W}^{(j)}_{\Z}(g,\lambda)$ \ $\cong$\ {\em Specht-Mod.}\  
${\cal S}^{\tau_j}_{\Z}$\ \ ($S_{n(\lambda)}$-module) \ \ 
\vs{3}

\noindent
\parbox[b]{1.1in}
{\ \ with 2-row

Young Diagram 

\ \ ($n(\lambda)$ boxes)}
\ \ \raise8pt\hbox{$\tau_j\,=\;$}
\parbox[t]{3in}
{$\byg & \cr  & \cr \eyg$  \raise8pt\hbox{ . . . }\qquad
$\underbrace{\mkern-30mu \byg & & \cr \cr \eyg\raise12pt\hbox{ . . . } 
\raise14pt\hbox{$\byg & & \cr  \eyg$} }_{j-1}$}. 

\vs{2}

\noindent
$\Longrightarrow$\ \  {\bf Complex} \ \ of $p$-modular $S_n$-reps.
\vs{2}

$\ldots\to {\cal S}^{\tau_{c_i}}_p\,\to\,\ldots\,\to\,
{\cal S}^{\tau_{2p-k}}_p\,\to\,{\cal S}^{\tau_{k}}_p\,\to\,
{\cal D}^{\tau_{k}}_p\,\to\,0$
\vs{2}

\noindent
{\bf Jordan H\"older (Composition)
Series:} \ \ 

$0\subset M_0\subset M_1\subset\ldots
\subset M_{m-1}\subset M_m={\cal S}^{\tau}_p$

All  $R_i=\raise2pt\hbox{$M_{i}$}\big/
\raise-2pt\hbox{$M_{i}$}\;\neq\;0$ and  irreducible

Write $R_0|R_1|\ldots|R_{m-1}|R_m$  and  ${\cal F}({\cal S}^{\tau}_p)=\{R_i\}$. 
\vs{2}

\noindent
{\bf Thm [James]} The ${\cal D}^{\tau}_p$ are irreducible. 

${\cal F}({\cal S}^{\tau_j}_p)\;=\;\{{\cal D}^{\tau_{j+2l}}_p\,:\;
l\geq 0\;\mbox{and}\;l\stackrel{p}{\subseteq} j+2l\}\;$
\vs{1.5}

\ \ All multiplicities one or zero. 

\vs{1}

\qquad\quad  ...\ \ \  not enough\ \ \  ...

\newpage 

\section*{Exactness via Modular Ordering:}

{\bf Partial Order} on ${\cal F}({\cal S}^{\tau}_p)$ (w/ mult=0,1) defined by 
\vs{2}

$R\leq R'$ \ \ iff \ \ \ $R_0|\ldots |R|\ldots |R'|\ldots |R_m$ \ \ $\forall$
JH-Series.
\vs{1}

\noindent
{\bf Thm [Kleshchev, Sheth, '99]} 

Determine\ \  $\leq$ \ \ for 2-row Young diagrams $\tau$ and $p\geq 3$:  
\vs{-2}

\noindent  
\ \ $({\cal F}({\cal S}^{\tau}_p),\geq)$ {\em ``inclusion structure''} 

\ w/ \  ${\cal D}^{\tau}_p$
maximal element. 
\vs{2}

\noindent
{\bf Cor [K]} Modular structure of maps in ${\cal C}(p,k)$
\vs{-2}

\picbox{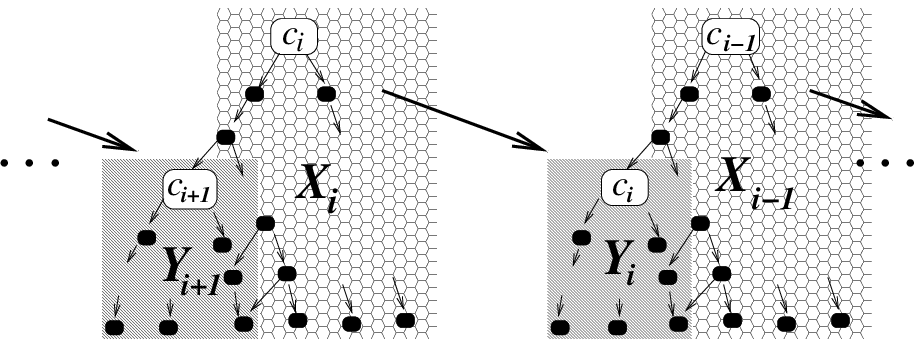}

with \ ${\cal F}(Y_i)\;=\;{\cal F}(X_i)$.\qquad\qquad
\vs{2}

\qquad\qquad\qquad\qquad\qquad $\Longrightarrow$\ \qquad \ EXACTNESS. $\ep$
\vs{3}

\noindent
{\bf Thm [K]} Image of Alexander Polynomial wrt. $\,t\,\mapsto\,\pm\zeta_p$: 
\vs{.5}

$\DP 
\dov\Delta_{\chi, p}^{\mp}(M)\;=\;\sum_{k=1}^{p-1}[k]_{\pm \zeta_p}
trace(\dov {\cal V}^{(k)}_p(C_{\Sigma}))\;\;\in\;\FF_p[\zeta_p]\;. 
$

\section*{Johnson-Morita Extension:}

{\bf Torelli Group:}\ \ \ \ $0\,\to\,{\cal I}_g\,\longrightarrow\,
\Gamma_g\,\longrightarrow\,Sp(2g,\Z)\,\to\,0$
\vs{3}

\noindent
{\bf Boundary Group:}\ \ \  ${\cal K}_g$  smallest normal containing

\
 
$D({\cal C}_k)$
\vs{-9}

\qquad\qquad\quad\epsfbox{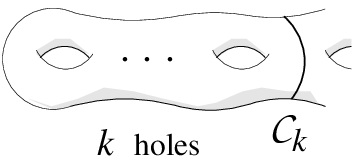}

\vs{2}

\noindent
{\bf Thm \small[Johnson, Morita]} $\,\exists\,$   Homomorphsims \ \ ($\scriptstyle H=H_1(\Sigma_g,\Z)$)
\vs{-1.5}

\noindent
\ $
\begin{array}{ccccccccc}
0 &\to&{\cal K}_g&\longrightarrow&
\Gamma_g&
\stackrel{\mbox{$\tilde k$}}{
-\!\!\!-\!\!\!\longrightarrow}
&\frac 12 \frac {\ext 3 \mbox{$H$}}{\mbox{$H$}}\rtimes Sp(2g,\Z)&\to&0\\
&&\Big \| &&\cup &&\cup&&\\
0&\to&{\cal K}_g&\longrightarrow&
{\cal I}_g&
\stackrel{\mbox{$\tau_2$}}{
-\!\!\!-\!\!\!\longrightarrow}&\frac {\ext 3 \mbox{$H$}}{\mbox{$H$}}
&\to&0
\end{array}
$
\vs{2}

\noindent
Let $g^*Jg=J$ for $g\in Sp(2g,\Z)$ and let $x\in \ext 3 H$.
\vs{2}

\noindent
Define   $\mu(x):\ext j H\to \ext {j-3} H$ by 
$\lz a,\mu(x)b\rz=\lz (Jx)\wedge a, b\rz$. 
\vs{0}

\noindent
{\bf Lem [K]} \ \ $\mu$ factors into $Sp(2g,\Z)$-covariant map
\vs{1}

\ $
\mu\,:\;\frac 12 \frac {\ext 3 \mbox{$H$}}{\omega\wedge \mbox{$H$}}
\;\longrightarrow\; 
Hom(\dov{\cal V}^{(k)}_p(\Sigma_g),\dov{\cal V}^{(k+3)}_p(\Sigma_g))
$
\vs{2}

\noindent
{\bf Thm [K]}\ For $\,0<k<p-3$ obtain \ \ 
$\raise3pt\hbox{$\Gamma_g$}\Big / \raise-3pt\hbox{${\cal K}_g$}$-module

$\dov {\cal U}^{(k)}_p(\Sigma_g)=
\dov{\cal V}^{(k)}_p(\Sigma_g)\stackrel{\mu}{\oplus}
\dov{\cal V}^{(k+3)}_p(\Sigma_g)$
\vs{-1}

\qquad\qquad\qquad by matrix \qquad$(x,g)\,\mapsto\,
\raise4pt\hbox{$\left [\begin{array}{cc}g&0\\
\mu(x)g & g
\end{array}\right]$}$

\newpage

\section*{Summary of TQFT Constructions:}

\bigskip

\bigskip

\begin{tabular}{rclc}
&${\cal V}^{FN}_\Z$ & Jacobian TQFT  &\\
&&                    Fully reducible over $\Z$ & \\
&&\ \ ({\em Permutation Modules})\vspace*{-.5cm} &\\
$\sl_2$-Lefschetz &\sdwa\vspace*{-.3cm}  & &  \\
Decomposition \vspace*{.2cm} & & & \\
&${\cal V}^{(k)}_\Z$ & Lefschetz Component, &  \\
&&                    irreducible  over $\Z$ \vspace*{-.05cm}&\\
&&\ \ ({\em Specht Modules  ${\cal S}_{\Z}^{(j)}$})\vspace*{-.5cm} &\\
$p$-Reduction &\sdwa\vspace*{.3cm}  & &  \\
&${\cal V}^{(k)}_p$ & Reducible over $\FF_p$ &  \\
&& with inner form.   \vspace*{-.04cm}&\\
&&\ \ ({\em $p$-Specht Modules  ${\cal S}_{p}^{(j)}$})\vspace*{-.4cm} &\\
Null space quotient &\sdwa\vspace*{.3cm}  & &  \\
&$\dov {\cal V}^{(k)}_p$ & Irreducible over $\FF_p$. &  \\
&&  Not $\Z$-liftable, but \vspace*{-.1cm} &\\
&& resolutions in ${\cal V}^{(j)}_p$'s   \vspace*{-.0cm}&\\
&&\ \ ({\em Simple $S_n$-modules  ${\cal D}_p^{\{j\}}$})\vspace*{-.7cm} &\\
Johnson-Morita &\sdwa\vspace*{-.3cm}  & &  \\
Extension \vspace*{.3cm}& & & \\
&$\dov {\cal U}^{(k)}_p$ & Indecomposable with &     \\
&&        two factors $\dov {\cal V}^{(k)}_p$ and $\dov {\cal V}^{(k+3)}_p$ \hspace*{-1cm}&\\
\end{tabular}
\bigskip

\newpage 

\section*{The Fibonacci TQFT:}

\noindent 
${\cal C}_5$ has two ``colors'' (simple objects) , namely \ 1\ and $\rho$\
\vs{2}

with \  \  $\rho\otimes\rho\;=\;1\oplus\rho$\  \ and \ $1$ \  = \ unit. 
\vs{1.5}

\noindent 
{\bf RT-Vectorspace :} \ \ over $\Q\,[\zeta_5]$ 
\vs{1.5}

${\cal V}^{RT}_{\zeta_5}(\Sigma_g)=Inv_{{\cal C}_5}({\bf F}^{\otimes g})$\ \ \ with \ \ \ 
${\bf F}\,=\,1\oplus(\rho\otimes\rho)
$.
\vs{1.5}

\noindent
{\bf Fork basis:} Example $g=10$  \ \ \ (1$\to$dotted; \ \ \ $\rho\to$ solid)
\vs{1}

\epsfbox{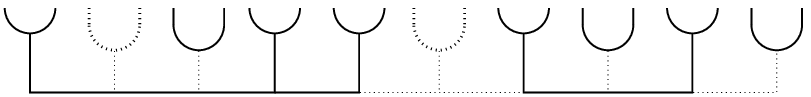}
\vs{1}

\noindent
{\bf Minimal Modules:}   

${\mathfrak S}_p(\Sigma)\;=\Z[\zeta_p]-\mbox{span of\ } \{
{\cal V}^{RT}_{\zeta_p}(M)\Omega\,:\;M:\emptyset\,\to\,\Sigma\}$
\vs{1}

\noindent
{\bf New Basis:}  \ \ Set $\;{\sf x}\,=\,\zeta_5-\zeta_5^{-1}\;$   \ \ and \ \ 
$D=-(\zeta_5^2+\zeta_5^{-2})\,$. 
\vs{-.5}

assign: \ \quad \ \raise-1cm\hbox{\epsfbox{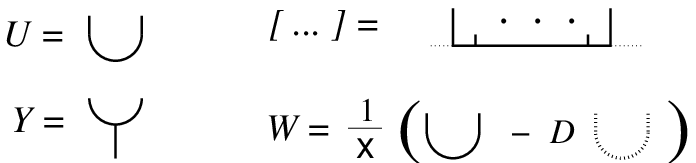}}

and rescale with \ \ \ $\DP {\sf x}^{\mbox{$[\frac {\#\;of\;Y's}2]$}}$\ \ \ obtain  new 
\vs{1}

$\Longrightarrow$ basis vectors $\lz\!\lz [Y,W,U,Y,Y][W][Y,U,Y][U]\rz\!\rz$ 
\vs{1}

{\em denote:}  
${\cal V}^I_{\zeta_5}(\Sigma)$ the $\Z[\zeta_5]$-span of these.
\vs{2}

\noindent
{\bf Thm [K]}   

$\exists_{\bf 1}$ \  \  {\em half-projective} TQFT ${\cal V}^I_{\zeta_5}$ with 
parameter ${\sf x}$  and 
\vs{1}

\noindent 
vectors spaces
\qquad\ \ 
${\cal V}^I_{\zeta_5}(\Sigma)\;=\;{\mathfrak S}_5(\Sigma)$
\vs{2}

{\em Proof:}  Explicit computation of $\Cob$-generators.

\newpage

\section*{Identification of $\FF_5$ Reductions:}

{\em Denote:} \ \ $\,{\cal V}^I_{0,5}$ \ \ the \ \ $\Z[\zeta_5]\to\FF_5$ \ \ 
reduction of $\,{\cal V}^I_{\zeta_5}$. 
\vs{3}

\noindent
{\em Recall:}  With $\,0\,<\,k=1\,<\,p-3\,=\,2$ have  construction
\vs{0.5}

\qquad
$\dov {\cal U}^{(1)}_5(\Sigma_g)\;=\; 
\dov {\cal V}^{(1)}_5(\Sigma_g)\,\stackrel{\mu}{\oplus}\, 
\dov {\cal V}^{(4)}_5(\Sigma_g)\,$
\vs{1.5}

\quad . . . \quad {\em guess what} \quad . . . 
\vs{2.5}

\noindent 
{\bf Thm [K]} \qquad $\dov {\cal U}^{(1)}_5\;\cong\;{\cal V}^I_{0,5}$\qquad
(over $\FF_5$).
\vs{2.5}

\noindent 
{\em Elements of Proof:}
\vs{2.5}

Middle dim  in $\ext * \lz a_i, b_i\rz$: \ \  $a_i\,\mapsto\,U$ and $b_i\,\mapsto\,W$ . 
\vs{2}

\noindent 
 $\dov {\cal V}^{(1)}\;\;\cong\;$ even \# of $Y$'s \qquad \quad 
 $\dov {\cal V}^{(4)}\;\;\cong\;$ odd\# of $Y$'s 
\vs{2}

\noindent 
Map $\raise-6pt\hbox{$\bsyg i \cr j\cr \esyg$}\;\leftrightarrow \;a_i\wedge b_i-a_j\wedge b_j$ 
\vs{-3}

\noindent 
to $Y$-pair in 
$i$-th and $j$-th position \quad \raise-0.2cm\hbox{\epsfbox{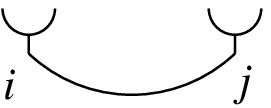}}
\vs{2}

\noindent 
Resolve via T.L.-Recoupling theory into fork-tree gadget. 

\vs{3}

\noindent
{\bf Dimensions: } \ {\em Fibonacci}  $\#$'s : \ \ $f_{n+1}\,=\,f_n+f_{n-1}$\  

(with \ \ $f_0=0$ , $f_1=1$, $f_2=1$, $f_3=2$,  . . . )
\vs{2}

\noindent 
For $g$ even\ \ $dim({\cal V}^I_{0,5}(\Sigma_g))\,=\,5^{\frac g2}f_{g-1}\;$ \ \ . . . 
\vs{2.5}

$\Longrightarrow$ \ \ From Resolution of $S_n$-{\bf\small [Ryba]}-reps. \   
 \ $\Longrightarrow$
\vs{2}

\noindent 
$f_{2r}\,=\,c(2r,r-1)-c(2r,r-3)+c(2r,r-6)-c(2r,r-8)+$

\qquad\  $+c(2r,r-11)-c(2r,r-13)+c(2r,r-16)-\ldots\,$
\vs{2}

{\em Catalan } \#'s: \  \ $c(n,j)={n\choose j}-{n\choose {j-1}}$.

\newpage

\section*{Applications of Fibonacci Identification:}

\begin{itemize}

\item 
$\DP 
2\dov\Delta_{\chi, 5}^{+}(M)\;=\;dim(ker(\chi))\,-\,
{\sf x}{\cal V}^{RT}_{\zeta_5}(M)\;+\;
{\cal O}({\sf x}^3)$.

follows from character formula for $\dov\Delta_{\chi, p}^{\pm}(M)$, and 
${\cal V}^{RT}(M)\,=\,{\sf x}\cdot trace ({\cal V}^{RT}(C_{\Sigma})$ and Thm. 

\noindent 
Note $\Rightarrow$ that if $\beta_1(M)>1$ have 

{\bf Thm[K]}\ \ \ \ $\raise4pt\hbox{$\dov\Delta_{\chi, 5}^{+}(M)$} \Big / \raise-4pt\hbox{${\sf  x}^3$}$ 
\ is\ \ {\em independent}  \ of \  $\chi$ !!

\ $\Rightarrow$\ Constraints on 3-manifold groups. 

\item If $\,\partial M\,=\,\emptyset\,$ then \ 
${\cal V}^{RT}_{\zeta_5}(M)=det(H_1(M))+{\cal O}({\sf x})$.

\ \ \ (= Invariant of $\dov {\cal V}^{(1)}$ and vacua lie in $\dov {\cal V}^{(1)}$)

\item The $\lambda_j^{O} {\,mod\,}  5 $ are of finite type \ \ (known {\bf \small [Kricker]}). 

{\em Different Proof:}{\small Pick element in augmentation ideal $\,u\,\in \,I{\cal I}_g\,$. 
Represented by ${\cal V}^I_{\zeta_5}$ \ in \ $\left[\begin{array}{cc}
0 & 0\\
\mu(u) & 0\end{array}\right]\,+\,{\sf x}\cdot End$. Thus $u^2\,\in\, {\sf x}\cdot End$
and $u^{2m}\,\in\, {\sf x}^m\cdot End$. 
Also $v\in {\cal K}_g\,$ is represented in ${\sf x}\cdot End$ so that 
$v^{m}\,\in\, {\sf x}^m\cdot End$. 

Hence $\lambda_{m-1}^{O}=0$ for splittings in $(I{\cal I}_g)^{2m}$ and 
$(I{\cal K}_g)^{m}$ $\Rightarrow$ finite type by {\bf \small[Garoufalidis, Levine]}. }

\item Have ${\cal V}^I_{\zeta_5}(D({\cal C}_k))\;=\;\I\,+\,{\sf x}\zeta_5^{n}\Pp_k$ 
with $\Pp_k^2=\Pp_k$. 

From $D(\psi({\cal C}_k))=\psi\circ D({\cal C}_k)\circ \psi^{-1}$ we find 

${\cal V}^{RT}_{\zeta_5}({D(\psi({\cal C}_k))})\,=\,
\I\,+\,{\sf x}\dov{\cal U}(\psi)\Pp_k\dov{\cal U}(\psi)^{-1}\,+\,{\cal O}({\sf x}^2)
$

{\bf\small [Murakami]}\ $\Rightarrow$\ ${\lambda}_{Casson}\,=\, 
2\lz\Omega, \dov{\cal U}(\psi)\Pp_k\dov{\cal U}(\psi)^{-1}\Omega\rz$
\end{itemize}
\newpage

\section*{Cut-Numbers \ \& \ Quantum-Orders:}
\vs{-1}

{\em Define:} \ \ $cut(M)$ \ for closed, connected \ $M$ \ by 

\begin{itemize}
\item Largest $k\in\N\cup\{0\}$ \ s/th\  $\exists$ 2-sided $\Sigma\subset M$

with \ \ $\pi_0(\Sigma)=k$ \ \ but \ \ $\pi_0(M-\Sigma)=1$\ . 

\item Largest $m\in\N\cup\{0\}$ \ s/th\  $\exists$ \ epimorphism 

$\;\pi_1(M)\,-\!\!\!\!\to\!\!\!\!\!\to\, F_m\,\cong\,\underbrace{\Z*\ldots*\Z}_m$\ \ 
onto free group. 
\end{itemize}
\vs{-1}

\noindent
{\bf Thm [?]}  \ \ Two definitions are equivalent. 
\vs{2}

{\em Note:}  \ \ $\beta_1(M)\,\geq\,cut(M)$ \ \ and \quad

\qquad \qquad \qquad $\beta_1(M)\neq 0 \,\Leftrightarrow \,cut(M)\neq 0$ \
\vs{1.5}

\nopagebreak
\noindent
{\bf Cor [K]}  \ \ Suppose ${\cal V}$ is half-projective w/ parameter ${\sf x}$
\vs{1}

\quad Then \ \ \ ${\cal V}(M)\;\in\;{\sf x}^{cut(M)}\cdot {\sf R}\;$.   
\vs{1.5}

\parbox[t]{1.3in}{
{\em Proof:} Present 

\noindent 
$\,M\,=\,A\circ B\,$ 

\ \ with 

\noindent 
$\,\mu(A,B)=cut(M)\,$}
\ \ \quad \raise-2cm\hbox{\epsfbox{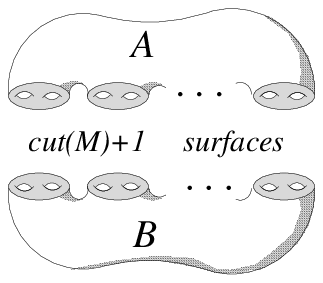}}
\vs{2}

\noindent 
{\em Definition:} Let ${\mathfrak o}_p(M)$ largest integer $k$ such that
\vs{1.5}

\qquad ${\cal V}^{RT}_{\zeta_p}(M)\;\in\;(\zeta_p-1)^k\Z[\zeta_p]$
\vs{2}

\noindent
{\bf Thm {\small [Gilmer, K]}}  \ \quad  \ $cut(M)\,\leq\,{\mathfrak o}_5(M)$. 
\vs{2}

\noindent
{\bf Thm {\footnotesize  [Melvin, Cochran]}}  
\ \quad  \ $\DP \frac {\beta_1(M)}3\,\leq\,\frac {{\mathfrak o}_p(M)}
{{\mathfrak o}_p(S^1\times S^2)}
$. 
\vs{2}

\newpage

\section*{Two Examples {\small [with  Gilmer]}: }

\noindent 
{\bf Example 1 (0-surgery along links):} 

\noindent 
\ Let $M({\cal L})$ be 3-mfld obtained via $0$-surgery along 
${\cal L}\subset S^3$ 
\vs{-2}

\noindent
{\bf Lem:} If $\cal L$ is boundary then $cut(M({\cal L}))=\,|{\cal L}|\,$. 
\vs{1.5}

\noindent 
{\em Denote:}${\cal B}={\cal C}_1\cup{\cal C}_2\cup{\cal C}_3$ 
{\em Borromean Link}. 

${\cal C}^*$ the {\em Whitehead Double} of a component $\cal C$. 

${\cal B}^*={\cal C}_1^*\cup{\cal C}_2\cup{\cal C}_3$,\ \ 
${\cal B}^{**}={\cal C}_1^*\cup{\cal C}_2^*\cup{\cal C}_3$, 
\vs{1.3}

\noindent 
Clearly, $\beta_1(M({\cal B}^{\ldots}))=3$. 
Note,  that $M({\cal B})=(S^1)^{\times 3}\,$ 

so that  $\pi_1(M({\cal B}))=\Z^3$
and hence  $cut(M({\cal B}))=1$. 
\vs{1.2}

\noindent
Moreover, ${\cal B}^{**}$ is boundary so that $cut(M({\cal B^{**}}))=3$. 
\vs{-1}

\noindent 
\parbox[t]{2.3in}
{\vs{1}

\ \ What is $cut(M({\cal B^{*}}))$ ?

\ \ Have two surfaces ${\cal R}_1$ 
and ${\cal R}_2$.

\ \ Compute ${\mathfrak o}_5(M({\cal B^{*}}))=2$

$\Rightarrow$ \ Answer \ \ $cut\;=\;2$ !

\smallskip

\noindent
{\em Note also:}\ \ ${\mathfrak o}_5(M({\cal B}))=1$
and

 \qquad   \qquad   \qquad  ${\mathfrak o}_5(M({\cal B^{**}}))=3$. 

}\quad 
\parbox[t]{2in}
{\vs{0}

\epsfbox{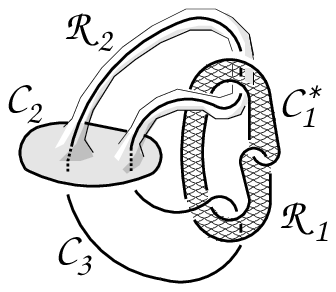}}
\vs{1.5}

\noindent 
{\bf Example 2 (Mapping Tori):}

For $j=1,\ldots,g-1$ \ let \  ${\cal A}_j\subset \Sigma_g$ be homologically
independent, mutually non-intersecting curves, and $D_j=D({\cal A}_j)$ the 
corresponding 
{\em Dehn twists}.

Let $M(n_1,\ldots,n_{g-1})\,\cong\,
\raise2pt\hbox{$\Sigma_g\times[0,1]$}\Big /
\raise-2pt\hbox{$(x,0)\sim (\psi(x),1)$}\,$ be the mapping torus 
with $\psi=D^{n_1}_1\circ \ldots\circ 
D^{n_{g-1}}_{g-1}\,$. 
\vs{1}

\noindent 
{\bf Prop:} If all $n_j\not\equiv 0 \, mod\, 5$ and $\,g\not\equiv
0 \, mod\, \,5$, 

then \qquad  $cut(M(n_1,\ldots,n_{g-1}))\;=\;g$
\newpage

\section*{Open Problems:}

\begin{itemize}

\item Find $\Z[\zeta_p]$-basis for ${\cal V}^{RT}_{\zeta_p}$  for $p\geq 7$. 

{\bf \footnotesize [Gilmer]}  Basis for $g=1$ and {\em all}  $p$. 

 \ \ (? w/ {\bf \footnotesize [Masbaum]} for $g=2$ ?)

\item Show they're half-projective \ w/\  ${\sf x}=(\zeta_p-\zeta^{-1}_p)^{\frac {p-3}2}$.

\quad $\Longrightarrow$ \ \ $cut(M)\,\leq\,\frac {{\mathfrak o}_p(M)}
{{\mathfrak o}_p(S^1\times S^2)}\,$. 

\item Characterize  $\,\DP \xi(M)=\frac {\beta_1(M)}{cut(M)}\,$! 
Is it bounded? 

When is $ \xi(M)\geq 3$? \ \ (Answers by 
{\bf \footnotesize [Sikora]})

\item Identify $\FF_p$-reduction ${\cal V}_{0,p}^I$ with 
$\dov {\cal V}^{(k)}_p$ and $\dov {\cal U}^{(k)}_p$. 

{\bf NoGo Thm:} For $p\geq 7$ cannote be linear relation.

{\em Proof:} by large $g$-asymptotics 

$dim(\dov {\cal V}^{(k)}_p(\Sigma_g))\sim  f_p^g$ and
$dim({\cal V}^{RT}_{\zeta_p}(\Sigma_g))\sim F_p^g$. 

Have $f_5=F_5$ but, e.g.,  $F_7=2f_7^2-7f_7+7$.

Suggests a form \ \ ${\cal V}_{0,7}^I\,\subset\, 
\FF_2\!\otimes\!\dov {\cal V}^{(\#)}_7
\!\otimes\!\dov {\cal V}^{(\#)}_7$. 

\item Is there TQFT for JM-extensions $\dov {\cal U}^{(k)}_p$ ? 
\ \ 
(By deformations of Jacobian-TQFT or  $\FF_2\ltimes\ext *\R^2$-theory?)  

\item What is structure of ${\cal V}_{j,p}^I$ for $j\geq 1$ ? 

$\Rightarrow$ \ \ Deeper structure of Mapping Class Group ... 

\item Use to find cohomology and characteristic 
classes of Mapping Class Group. 

\item For $\beta_1(M)>0$ is ${\cal V}_{1,p}^I$ given by 
$\dov\Delta_{\chi}(M)$.\ \ How?! 

Relation to Frohman-Nicas-$PSU(n)$-Knot Invariants? 

\item Which parts of $\dov\Delta_{\chi,p}(M)$ independent of $\chi$?
 
Interpretation by $\FF_p$-covering, signature defects, etc.? 

\item What families of algebraic groups can be excluded from this
as 3-manifold groups ?  

\item See ${\cal V}_{j,p}^I$  as ``higher torsion'' - invariants ?  

\quad . . . as ``$\beta_1>0$''-extensions of $\lambda_j^O$? 

\item Relations to 3D gauge theory (Seiberg Witten, Donaldson)

\item Reexpress formulae for $\lambda_{Casson}\;mod\;5\;$ 
for Heegaard splits with $g\geq 2$
in $Sp(2g,\Z)$-invariant-theoretical terms. 

Compare to Morita's formulae!

\ \ For $g=2$ know $\Pp_1$ explicitly \ 

\ \ $\Rightarrow$\ precisely reproduce 
Morita formula. 
\vspace*{-.1cm}

\ \ \  \ \ \ {\small (looks a-priori different for $g>2$)}

\item Same for general $\;p\;$. 

\item Heegaard formulae also for $\lambda_j^O\;mod\;p$ for $j>1$ and 
general $p$. 

\end{itemize}

\end{document}